\documentstyle[11pt,pb-diagram,a4,twoside]{article}
\newtheorem{theorem}{Theorem}[section]
\newtheorem{lemma}[theorem]{Lemma}

\newtheorem{proposition}[theorem]{Proposition}
\newtheorem{satzdef}[theorem]{Theorem/Definition}

\newcommand{\rk}{{\rm rk}}
\newcommand{\Hom}{{\rm Hom}}

\newcommand{\length}{{\rm length}}
\newcommand{\Pic}{{\rm Pic}}

\newcommand{\Ext}{{\rm Ext}}

\newcommand{\Spec}{{\rm Spec}}

\newcommand{\U}{{\rm U}}
\newcommand{\SU}{{\rm SU}}
\newcommand{\NS}{{\rm NS}}

\newcommand{\gr}{{\rm gr}}

\newcommand{\triv}{{\rm triv}}
\newcommand{\supp}{{\rm supp}}

\newcommand{\Ecal}{{\cal E}}

\newcommand{\Gcal}{{\cal G}}

\newcommand{\Jcal}{{\cal J}}
\newcommand{\Lcal}{{\cal L}}

\newcommand{\Ocal}{{\cal O}}

\newcommand{\pdop}{I \!\! P}
\newcommand{\zdop}{Z \!\!\! Z}

\newcommand{\dual}{^\lor}
\newcommand{\ddual}{^{\lor \lor}}

\newcommand{\rar}{\rightarrow}
\newcommand{\Rar}{\Rightarrow}
\newcommand{\rarpa}[1]{\stackrel{#1}{\rightarrow}}

\newcommand{\larpa}[1]{\stackrel{#1}{\leftarrow}}

\newcommand{\remark}{{\bf Remark:  }}

\newcommand{\proof}{{\bf Proof:  }}
\newcommand{\definition}{{\bf Definition:  }}
\newcommand{\case}[1]{{\bf Case #1:  }}
\newcommand{\step}[1]{{\bf Step #1: \hspace{0.2em} }}
\newcommand{\qed}{{ \hfill $\Box$}}
\newcommand{\tiptop}{{\vspace{0.5em} }}

\begin{document}

\title{Duality Construction of Moduli Spaces}
\author{Georg Hein}
% \date{8.~April 1997}
\maketitle

%\tableofcontents

\section*{Introduction}
In \S 1 of Faltings' article \cite{Fal} a ``GIT-free'' construction
is given for the moduli spaces of vector bundles on curves using
generalized theta functions.
Incidentally, this construction is implicitly described
in Le Potier's article \cite{LP2}.
The aim of this paper is to generalize the
{\em duality construction}
to projective surfaces.

For a rank two vector bundle $E$ on the projective plane $\pdop^2$,
the divisor $D_E$ of its jumping lines is a certain generalization
of the Chow divisor of a projective scheme.
We give a generalization of this divisor for coherent sheaves on
surfaces.
Using this duality we construct the moduli space of coherent sheaves
on a surface that does not use Mumford's geometric invariant theory.
Furthermore, we obtain a finite morphism from this moduli space to
a linear system, which generalizes the divisors of jumping lines.

Applying this construction to curves, we get exactly Faltings'
construction.
The moduli space we construct here can also be obtained by using GIT.
This construction is carried out in \S8.2 of the
book \cite{HL} of Huybrechts and Lehn.
Le Potier obtained this moduli space in \cite{LP1} for surfaces
with ``many lines'' (see \S \ref{SBAR} for an exact definition).
However, it is the modest hope of the author that the construction
presented here provides new insight into the geometry of moduli spaces.

First we outline this concept, which generalizes the famous
{\em strange duality} to moduli of coherent sheaves on surfaces.
To do so we define duality between schemes in part \ref{SDUAL},
giving three examples of ``natural dualities''.
In section \ref{SDCON} the duality construction is given.
In order to avoid a too technic presentation of the construction itself,
we defer the proofs to the following section.
The last section is dedicated to the Barth morphism.

In order to simplify the discussion we restrict ourselves to
moduli spaces of sheaves of rank two with trivial determinant.
The interested reader will be able to extend this
to arbitrary rank and determinant.

The author is thankful to his thesis advisor,
H.~Kurke, for many fruitful discussions.

\section{Duality of schemes}\label{SDUAL}
\subsection{Definitions}
Let $(X,\Ocal_X(D_X))$ and
$(Y, \Ocal_Y(D_Y))$ be two schemes with line bundles.
A duality between these two pairs is given by a nontrivial
section $s \in H^0(\Ocal_X(D_X)) \otimes H^0(\Ocal_Y(D_Y))$.
We will identify $s$ with its vanishing divisor
$D=V(s) \subset X \times Y$.

\vspace{1em}
\hspace{12em}
\setlength{\unitlength}{0.0012cm}
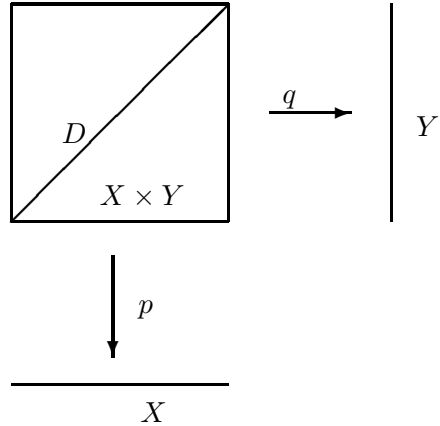
\begin{picture}(4512,4548)(2389,-6097)
\thicklines
\put(4801,-1561){\line( 0,-1){2400}}
\put(4801,-3961){\line(-1, 0){2400}}
\put(2401,-3961){\line( 0, 1){2400}}
\put(2401,-1561){\line( 1, 0){2400}}
\put(4801,-1561){\line(-1, 0){2400}}
\put(2401,-3961){\line( 1, 1){2400}}
\put(6601,-1561){\line( 0,-1){2400}}
\put(2401,-5761){\line( 1, 0){2400}}
\put(5251,-2761){\vector( 1, 0){900}}
\put(3526,-4336){\vector( 0,-1){1125}}
\put(3100,-2986){\makebox(0,0)[cc]{$D$}}
\put(5476,-2611){\makebox(0,0)[cc]{$q$}}
\put(6999,-2911){\makebox(0,0)[cc]{$Y$}}
\put(3976,-6061){\makebox(0,0)[cc]{$X$}}
\put(3888,-4936){\makebox(0,0)[cc]{$p$}}
\put(3851,-3686){\makebox(0,0)[cc]{$X \times Y$}}
\end{picture}\\
We obtain a rational morphism
$$\begin{array}{cccc}
s_X: & X & - - - \rar & |D_Y| \\
& x & \mapsto & q(D \cap p^{-1}(x)) \, . \\
\end{array}$$
The base locus of this morphism $s_X$ consists of all points $x$ of
$X$ such that the vertical component $x \times Y$ is contained in $D$.
This motivates the following

\tiptop
\definition
The duality $D$ between $(X,D_X)$ and $(Y,D_Y)$ is called\\
{\em generated,} if $s_X$ is a morphism;\\
{\em generated ample,} if $s_X$ is a finite morphism;\\
{\em very ample,} if $s_X$ is an embedding.\\

\subsection{Examples}
The first example demonstrates that the above definitions
are something with which we are familiar.

\tiptop
{\bf The duality of a linear system}
Let $X$ be a given scheme with an effective divisor $D_X$,
and let $Y \subset |D_X|$ be a linear system.
We take $D$ to be the incident divisor, i.e.
$D= \{ (x, H) \, | x \in H \}$.
Then the notions for $D$ given in the above definition
correspond to those for the linear system. 

\tiptop
{\bf Strange duality}
Let $C$ be a smooth projective curve of genus $g$
over the complex numbers.
We fix two positive integers $m$ and $n$
and a theta characteristic $A$,
i.e. $A \in \Pic^{g-1}(C)$ and $A^{\otimes 2} \cong \omega_C$.
We consider the following moduli schemes.
$$\begin{array}{ccl}
X & = & \U_C(n,n(g-1)) \\
 & = & \left\{ E \, \left| \,
\begin{array}{l}
E \mbox{ semistable } C \mbox{-vector bundles with } \\
\rk(E)=n \mbox{ and } \deg(E)=n(g-1)  \\
\end{array}
\right. \right\} \\
D_X & = & \{ E \in X \, | h^0(E)=h^1(E)>0 \} \\
Y & = & \SU_C(m) \\
 & = & \left\{ F \, \left| \, 
\begin{array}{l}
F \mbox{ semistable } C \mbox{-bundles with } \\
\rk(E)=m \mbox{ and } \det(E) \cong \Ocal_C  \\
\end{array}
\right. \right\} \\
D_Y & = & \{ F \in Y \, | h^0(F \otimes A)=h^1(F \otimes A)>0 \} \\
D & = & \{ (E,F) \, | \, h^0(E\otimes F) = h^1(E \otimes F) >0 \} \\
\end{array}$$
The line bundle $\Ocal_{X \times Y}(D)$ is isomorphic to
$p^*\Ocal_X(D_X)^{\otimes m} \otimes q^*\Ocal_Y(D_Y)^{\otimes n}$.
By the duality we obtain a linear map from
$H^0(X,\Ocal_X(D_X)^{\otimes m})\dual
\rarpa{s} H^0(Y,q^*\Ocal_Y(D_Y)^{\otimes n})$.
According to the Verlinde formula, both spaces have the same dimension.
The natural conjecture that $s$ is an isomorphism
is called the {\em Strange Duality Conjecture}.
For more details on this topic, see Beauville's survey article
\cite{Bea}.

\tiptop
{\bf Duality between moduli spaces on polarized surfaces}
The next example is our main example,
the notions of which will be used for the remainder of the article.
We describe the moduli spaces here only set theoretically.
Their construction uses the duality,
thus giving a rough idea of the construction
which is the object of the following section.
The concept of semistability used here is the Mumford (slope)
semistability.

Let $(S,\Ocal_S(1) = \Ocal_S(H))$ be a projective polarized surface.
Fix a class $c_2 \in H^4(S,\zdop)$.
We will consider a duality between the following two coarse moduli
spaces:
$$\begin{array}{rcl}
X & = & M_S(2,0,c_2) \\
& = & \left\{ E \left| \,
\begin{array}{l}
E \mbox{ semistable torsion free sheaf on } S, \mbox{ with} \\
\rk(E)=2 \quad \det(E) \cong \Ocal_S \quad c_2(E)=c_2 \\
\end{array}
\right. \right\} \\
\\
Y & = & M_{|H|}(2,\omega_{|H|}) \\
& = & \left\{ F \left| \,
\begin{array}{l}
F \mbox{ semistable torsion sheaf on } S, \mbox{ with} \\
Z=\supp(F) \in |H| \quad \rk_Z(F)=2 \quad \det_Z(F) \cong \omega_Z \\
\end{array}
\right. \right\} \\
\end{array}$$

The duality will be given by the $X \times Y$ divisor
$$D= \{  (E,F) | H^*(E \otimes F) \not= 0 \}\, .$$

\section{The duality construction}\label{SDCON}
Using the notation introduced in the last example,
we give a construction of the coarse moduli scheme
$X$ using the duality morphism $s$.
$X$ will be obtained together with a polarization
and the {\em Barth morphism}, which will be finite by construction. 
The steps for this construction are listed below,
and proofs for all pertinent theorems are provided in the next section.

\begin{description}
\item [Boundedness of $X$]

There exists a projective scheme $Q$ and a torsion free sheaf
$\Ecal$ on $Q \times S$
flat over $Q$ which (over)parameterizes the moduli problem.
More precisely denote by $p$, $q$ the two projections
$$Q \larpa{p} Q \times S \rarpa{q} S \, .$$
For any sheaf $E$ of $X$ let $Q_E$ be the subscheme
$$Q_E = \{ q\in Q \, | \, \Ecal_q \cong E \}$$
of $Q$. All $Q_E$ are required to be connected and nonempty.
Since this is the same starting point like in the GIT construction,
it is obvious that $Q$ can be taken to be a suitable
Quot scheme (see \cite{Gro}).

\item [Elements of $Y$ give sections in a $Q$-line bundle $\Lcal$]
(see \S \ref{P1})

We will show that there exists a $Q$-line bundle $\Lcal$
and a global section $s_F \in H^0(Q, \Lcal)$, for any $F \in Y$.
The vanishing locus of $s_F$ is given by
$$V(s_F) = \{ q \in Q \, | \, H^*(S, \Ecal_q \otimes F) \not= 0 \} \,
.$$

\item [Base points correspond to unstable objects] (see \S \ref{P2})

The base locus $B(\Lcal)$ with respect to the sections given by $Y$
is the scheme theoretic intersection of all $V(s_F)$ for all $F \in Y$.
Since $Q$ is noetherian we can write
$$B(\Lcal) = \bigcap_{i=0}^N V(s_{F_i}) \, .$$
It will be shown that $B(\Lcal)$ consists exactly of those points $q \in
Q$
for which the sheaf $\Ecal_q$ is not semistable.

\item [Properness of $X$]
(see \S \ref{P3})

We have to show that semistable limits of semistable families exist.

\item [The line bundle $\Lcal$ is $X$-positive]
(see \S \ref{P4})

It will be shown that the degree of $\Lcal$ on a curve $C$
parameterizing
semistable objects is zero only if
the curve parameterizes Jordan-H\"older equivalent sheaves.

\item [The duality construction]
The rational morphism $s=(s_{F_0}: \ldots:s_{F_N})$ from $Q$ to
$\pdop^N$
leads  to a morphism
$\bar Q \rarpa{\varphi} \pdop^N$ after a blow up of $Q$. 
We consider the following diagram  
$$\begin{diagram}
\node{\bar Q} \arrow{s,r}{\pi} \arrow{se,t}{\varphi}
\arrow{e,t}{\varphi_0}
\node{X} \arrow{s,r}{\varphi_1} \\
\node{Q} \arrow{e,t,..}{s} 
\node{\pdop^N}
\end{diagram} \quad.$$
Here $\varphi = \varphi_1 \circ \varphi_0$
is the Stein factorization of $\varphi$.
Hence the Barth morphism $\varphi_1$ is finite.
By the above, any point of $X$ corresponds to exactly one
Jordan-H\"older
equivalence class of semistable bundles.
\end{description}

\section{Details and proofs}
\subsection{The $Q$-line bundle $\Lcal$ and its invariant
sections}\label{P1}

$\Lcal$ is defined to be the determinant bundle
$\Lcal = det(p_!(\Ecal \otimes q^*F)^{-1}$.
The definition does not depend on the choice of $F \in Y$,
because these elements coincide in the Grothendieck group $K_0(S)$.

\tiptop
{\bf Framed elements of $Y$ give sections in $\Gamma(\Lcal)$}\\
Let $F$ be a semistable element in $Y$ with support $Z \in |H|$.
Then there exists a short exact sequence
$$0 \rar F \rar \Ocal_Z(M)^{\oplus 3}
\rarpa{\alpha} \omega_Z\dual(3M) \rar 0 \, .$$
We remark that $M>>0$ can be chosen for all $F \in Y$.
Define $s_F$ to be the section
$\det(Rp_*(\Ecal \otimes q^* \alpha))$.
It is clear from the construction
that the vanishing divisor $V(s_F)$
is supported on those $q\in Q$, for which
$H^*(\Ecal_q \otimes q^*F)$ is not zero. 

\remark By abuse of notation we simply write $s_F$
and do not explicitly refer to the framing.

\tiptop
{\bf Global sections in $\Lcal^{\otimes k}$}

Let $\tilde F$ be a rank two vector bundle on a curve
$\tilde Z$ from the linear system $|kH|$.
We require the determinant of $\tilde G$
to be isomorphic to $\omega_{\tilde Z}$.
Using adjunction to express $\omega_{\tilde F}$
and $\omega_F$,
the following computation in the Grothendieck group
$K_0(S)$ shows
that $[\tilde F] =k[F]$ for any $F \in Y$:
$$\begin{array}{ccl}
[\tilde F] & = & ([\Ocal_S]-[\Ocal_S(-kH)])([\Ocal_S]
+ [K_S(kH)]) \\
\\
& = & ([\Ocal_S] -[\Ocal_S(-H)])
\left(\sum\limits_{i=0}^{k-1}[\Ocal_S(-iH)])
([\Ocal_S] + [K_S(kH)] \right) \\
\\
& = & ([\Ocal_S] -[\Ocal_S(-H)]) k([\Ocal_S] +
[K_S(H)]) \, = \, k[F] \, \, .\\
\end{array}$$
Hence $\tilde F$ (together with a framing)
defines a global section in $\Lcal^{\otimes k}$.

\subsection{Semistability}\label{P2}
Assume that $H$ is big enough,
which means that $\Ocal_S(H)$ is globally generated,
and the following two conditions hold:

(i) $H^2>4c_2$, and

(ii) The positive generator $a$ of the
$\zdop$ ideal $\{ D.H \, | \, D \in \NS(S) \}$
satisfies $a>c_2$.

Under these assumptions we have the following
\begin{satzdef}\label{SEMISTABLE}
For a torsion free rank two $S$-sheaf $E$ with
$\det(E) \cong \Ocal_S$ and $c_2(E)=c_2$ the following
four conditions are equivalent:
\begin{enumerate}
\item There exists an $F \in Y$ such that $H^*(S,E\otimes F)=0$;
\item For all rank one subsheaves $M \subset E$,
the inequality $c_1(M).H \leq 0$ holds;
\item The restriction of $E$ to a general divisor $Z \in |H|$ is
semistable,
i.e. all $Z$-line bundles contained in $E_Z$
have nonpositive degree;
\item For $Z \in |H|$ general, there exists a $Z$-line bundle
$A$ such that $H^*(Z, L \otimes E|_Z)=0$.
\end{enumerate}
If one of these conditions is satisfied we call $E$ a semistable
$S$-sheaf.
\end{satzdef}

\proof {\bf (1) $\Rar$ (2) } Suppose there exists an
$M \subset E$ such that $c_1(M).H>0$.
Then $M$ restricted to $Z = \supp(F)$ is of positive degree.
Therefore the Euler characteristic of $M \otimes F$ is positive.
Since the sheaf is one-dimensional there are global sections.
Hence there are global sections in
$H^0(E \otimes F)$,
which contradicts the assumption of (1).

{\bf (2) $\Rar$ (3) }
This is a restriction theorem that follows from
Bogomolov's inequality (see \cite{Bog}). 
For a complete proof of this implication
see \cite{HL} theorem 7.3.5.

{\bf (3) $\Rar$ (4)} This result goes back to
Raynaud (\cite{Ray}).
For a shorter proof see \cite{He1}.

{\bf (4) $\Rar$ (1) }
Denote the genus of $Z$ by $g$.
It follows by Riemann-Roch that $A \in \Pic^{g-1}(Z)$.
The condition $H^*(E|_Z \otimes A) =0$ is satisfied on a
nonempty open subset of the Jacobian $\Pic^{g-1}(Z)$.
Hence the condition
$H^*(E|_Z \otimes \omega_Z \otimes A^{-1})=0$
is again open and not empty.
Consequently we can choose $F$ to be a direct sum
$A \oplus (\omega_Z\otimes A^{-1})$. \qed

\subsection{Properness}\label{P3}
Although the properness of the moduli functor $X$ is well known,
a new proof is given here which is shorter
than Langton's original proof in \cite{Lan}.
The main idea to get ``more and more'' semistable extensions
of a generic semistable family
by elementary transformations comes from Langton's proof.
However, using the invariant functions,
we can control the maximal number of elementary transformations
required.
Therefore the proof fits into the concept
of the duality construction.

\begin{theorem}
(\cite{Lan})
Let $R$ be a discrete valuation ring with $\Spec(R)=\{ 0, \eta \}$.
Let $\Ecal_\eta$ be a semistable torsion free sheaf on $\eta \times S$.
Then there exists an extension $\Ecal_R$ of
$\Ecal_\eta$
which is semistable in the special fiber as well.
\end{theorem}
\proof
We consider the following morphisms:
$$ \Spec(R) \stackrel{p}{\leftarrow}
\Spec(R) \times S \stackrel{q}{\rar} S \quad .$$
Since the Quot scheme is projective,
there exist torsion free extensions of $\Ecal_\eta$.
For an extension $\Ecal$,
we define its badness\footnote{We use the word badness
because $b$ measures
how far $\Ecal$ is from being a semistable extension.
So badness zero implies semistability.}
$b(\Ecal,F)$ with respect to an $F \in Y$ as
$$b(\Ecal,F)= \left\{ 
\begin{array}{ll}
\infty & \mbox{ if }
\supp R^1p_*(\Ecal \otimes q^*F)=\Spec(R) \, , \\
\length(R^1p_*(\Ecal \otimes q^*F)) & \mbox{ otherwise.} \\
\end{array}
\right.$$
The absolute badness $b(\Ecal)$ of $\Ecal$
is defined to be the minimum of all these numbers:
$$ b(\Ecal) = \min_{F \in Y} \{ b(\Ecal,F) \}  \, .$$
Since $\Ecal_\eta$ is semistable the badness
$b(\Ecal)$ has to be finite.
We suppose that $\Ecal$ is an extension with
minimal possible badness.
If the badness is zero, the special fiber
$\Ecal_0$ is semistable by \ref{SEMISTABLE}.
Hence we may assume that $b(\Ecal)>0$.
Since $\Ecal_0$ is not semistable,
there is a surjection
$\Ecal_0 \rar L \otimes \Jcal_Z$ with $L.H< 0$,
$\Jcal_Z$ being the
ideal sheaf of a codimension two subscheme of $S$.
We choose an element $F \in Y$
subject to the following three open conditions:

\tiptop
(i) $H^0(L \otimes F) =0$,

(ii) $b(\Ecal,F) = b(\Ecal)$,

(iii) $\supp(F) \cap Z = \emptyset $.

\tiptop
Define the elementary transformation $\Ecal'$ of
$\Ecal$ by the exact sequence
$$0 \rar \Ecal' \rar \Ecal \rar L \otimes \Jcal_Z \rar 0 \, .$$
Applying the functor $p_*(- \otimes q^*F)$ to that sequence,
we obtain
$$\begin{array}{cccc}
$$p_*(L \otimes q^*F) & \rar R^1p_*(\Ecal' \otimes q^*F) \rar 
R^1p_*(\Ecal \otimes q^*F) \rar 
& R^1p_*(L \otimes q^*F) & \rar 0\, .\\
|| && \parallel \!\!\!\!\!\! - \\
0 && 0 \\
\end{array}$$
This contradicts the minimality assumption
on the badness of $\Ecal$ by
the very definition of this number. \qed

\subsection{$X$-positivity of the line bundle $\Lcal$}\label{P4}
We have to consider the equivalence classes
of semistable sheaves parameterized
by our moduli space $X$ using the following equivalence relation.

\tiptop
{\bf Definition (trivially connected equivalence)}
Two semistable $X$-sheaves $E$ and $E'$ on $S$
are called trivially connected if
there exists a connected projective curve $B$
and a family $\Ecal$ on $B \times S$ such that\\
- the determinant line bundle $\Lcal_B$ on $B$ is trivial and\\
- there are points $b$ and $b'$ in $B$ with
$E \cong \Ecal_b$ and $E' \cong \Ecal_{b'}$.

\tiptop
There is a second equivalence relation
that reflects the geometry of the sheaves.
We start with some preparations.
If $\tau$ is a coherent sheaf of dimension zero,
then we define its trivialisation $\triv(\tau)$ by
$$ \triv(\tau) :=
\bigoplus_{P \in X} k(P)^{\oplus \length_P(\tau)} \, .$$
For a torsion free sheaf $G$,
let $G\ddual$ be its double dual and $\tau(G)$
be the cokernel of the injection
$G \hookrightarrow G\ddual$.
Define by $$\triv(G) = G\ddual \oplus \triv(\tau(G))$$
the trivialisation\footnote{This definition is good enough
for our purposes.
However it should be replaced by
$\triv(G) = G\ddual \ominus \triv(\tau(G))$.} of $G$.

\tiptop
Define the graded object of a stable sheaf
$E$ to be $E$ itself: $\gr_H(E) =E$.

If $E$ is a semistable but not stable sheaf,
then there exists a short exact sequence
$0 \rar A' \rar E \rar A'' \rar 0$
with $A'$ a saturated subsheaf of $E$ and $c_1(A').H=0$.
In this case we define the graduated object
$\gr_H(E)$ of $E$ to be the direct sum $A' \oplus A''$.

\tiptop
{\bf Definition (Jordan-H\"older equivalence)}
Two $X$-sheaves $E$ and $E'$ on $S$ are called
Jordan-H\"older equivalent if and only
if $\triv(\gr_H(E)) \cong \triv(\gr_H(E'))$.

\tiptop
This definition implies, in particular,
that the equivalence class of a stable vector
bundle consists of one element up to isomorphism.
In the course of the next result,
it will be shown that Jordan-H\"older equivalence
coincides with trivially connected equivalence.

\begin{theorem}\label{POSITIVE}
Let $B$ be a smooth projective connected curve and
$\Ecal_B$ be a family of sheaves on $B \times S$.
Assume that the sheaf parameterized by the
generic point of $B$ is semistable.
Then the $B$-line bundle $\Lcal_B$  is trivial or ample.
If $\Lcal_B$ is trivial,
then all $S$-sheaves parameterized by $B$ are Jordan-H\"older
equivalent.
\end{theorem}

First we note that there are nontrivial sections in $\Lcal_B$,
because there are semistable objects parameterized by points of $B$.
Therefore this line bundle is either trivial or ample.
Consequently the proof reduces to showing that
$\Ocal_B \cong \Lcal_B$ (by definition,
the trivially connected equivalence) induces the Jordan-H\"older
equivalence.
To prove the above theorem some preparations will be needed.
We retain the above notations.

\begin{theorem}\label{POSCURVE}
(\cite{Fal} theorem I.4)
Let $C$ be smooth projective curve and $\Ecal$
a vector bundle on $B \times C$ with
$deg_C(\Ecal_b) = 0$ for all $b\in B$.
Denote the projections of $B \times C$
to the components by $p$ and $q$.
Suppose there exists a $C$-line bundle $M$ such that
$R^*p_*(\Ecal \otimes q^*M) =0$ holds.
Then the $C$-objects parameterized by $B$ are all $S$-equivalent. 
\end{theorem}
\proof We proceed in steps.

\step{1} There exists a $B$-bundle $G$ such that
$\Ecal_P \cong G$ for all points $P \in C$.

Let $P$ and $Q$ be two arbitrary points of $C$.
The set of all line bundles $M$ in $\Pic^{g_C-1}(C)$
with $R^*p_*(\Ecal \otimes q^*M) =0$ is open.
Hence there is a line bundle $\tilde M \in \Pic^g(C)$
such that $\tilde M(-P)$ and $\tilde M(-Q)$ are in this open set.
From the exact sequence
$0 \rar \tilde M(-P) \rar \tilde M \rar k(P) \rar 0$ we obtain 
$$p_*(\Ecal \otimes q^*M) \cong
p_*(\Ecal \otimes q^*k(P)) \cong
p_*(\Ecal|_{  B \times \{ P \} })
\cong \Ecal|_{  B\times \{ P \} }\, .$$
Analogously, $p_*(\Ecal \otimes q^*M)
\cong \Ecal|_{ B\times \{ Q \} }$
which proves the assertion of the first step.

\step{2} Set $\Gcal = p^*G$.

There are three distinct cases to be considered.

\case{1} $G$ is stable

Since $G$ is simple it follows that
$N = q_*(\Gcal\dual \otimes \Ecal)$ is a $C$-line bundle.
But $\Gcal \otimes q^*N$ is isomorphic to $\Ecal$,
therefore all objects parameterized by $B$
are isomorphic to $N \oplus N$.

\case{2} $G$ is semistable but not stable

After a twist with a line bundle we may assume $G$
to be of degree zero.
By theorem \ref{SEMISTABLE} there exists a
$B$-line bundle $A$ such that
$G \otimes B$ has no cohomology.
This implies $R^*q_*(\Ecal \otimes p^*A) =0$,
and, as in the first step,
all $C$-objects parameterized by $B$ are isomorphic.

\case{3} $G$ is not semistable

Let $0 \subset A \subset G$
be the Harder-Narasimhan filtration of $G$,
i.e. $A$ is the subline bundle of $G$ of maximal degree.
We denote the quotient $G/A$ by $A'$.
By the uniqueness of $A$,
$N = q_*(p^*A\dual \otimes \Ecal)$ is a line bundle on $C$.
We find that $p^*A \otimes q^*N$ is a subbundle of
$\Ecal$ with cokernel isomorphic to
$p^*A' \otimes q^*N'$ for a $C$-line bundle $N'$.
Using the short exact sequence
$$ 0 \rar p^*A \otimes q^*(N \otimes M) \rar
\Ecal \otimes q^*M \rar p^*A' \otimes q^*(N' \otimes M) \rar 0$$
to compute the degree of $R^*p_*(\Ecal \otimes q^*M)$, we have
$$0 = \deg(A)\deg(N)+\deg(A')\deg(N') \, .$$
Since $\Ecal$ is a family of degree zero sheaves on $C$,
it follows that $\deg(N') = - \deg(N)$.
Hence we obtain the equality
$$0 = (\deg(A) - \deg(A'))\deg(N) \, .$$
By assumption the first factor is strictly positive,
thus $\deg(N) =\deg(N')=0$.
But this means that all objects parameterized by $B$ are
extensions of two line bundles of the same degree. \qed

\begin{lemma}\label{BOUNDLINE}
The set of all line bundles $L$ on $S$ such that $L.H=0$ and for
which there exists a nontrivial homomorphism in
$\Hom(E,L)$ for some $E \in X$ is bounded.
Subsequently, these line bundles can be parameterized by
an noetherian scheme.
\end{lemma}
\proof
It is enough to show that the set of Hilbert polynomials of
these line bundles $L$ is finite.
For any such line bundle $L$, there is an exact sequence
$$0 \rar L^{-1} \otimes \Jcal_{Z_1} \rar E \rar
L \otimes \Jcal_{Z_2} \rar 0 \, .$$
Since $E$ is semistable, $\Jcal_{Z_2}$ is the ideal sheaf of
some zero dimensional scheme.
Using the above sequence, the second Chern class can be computed,
and indeed, 
$c_2=-L^2+\length(Z_1) +\length(Z_2)$.
Hence we conclude by the Hodge index theorem that
$L^2$ is in the interval $[-c_2,0]$
and that $(H.(K_S \pm L))^2 \leq H^2 (K_S \pm L)^2$,
which gives lower and upper bounds for $K_S.L$.
The Hilbert polynomial of $L$ with respect to $H$ is determined
completely by
the numbers $L^2$, $L.K_S$ and $L.H$. \qed

\begin{lemma}\label{BIGK}
There exists a positive number $k$ such that, for all
$X$ sheaves  $E$ and $E'$ on $S$ and
all line bundles  on $S$ with $H.L=0$ and $\Hom(E,L) \not=0$,
the groups $\Ext^1(E, L(-kH))$ and $\Ext^1(E,E'(-kH)\ddual)$ vanish.
\end{lemma}
\proof
By Serre duality $\Ext^1(E,L(-kH)) \cong H^1(E(K_S+kH) \otimes
L^{-1})\dual$.
For any pair $(E,L)$ there exists a number
$k$ such that the cohomology group vanishes.
By lemma \ref{BOUNDLINE} the set of all these pairs is bounded.
Hence there exists a global $k$.

The same argument shows the vanishing of
$\Ext^1(E,E'(-kH)\ddual)$ for a given $k$. \qed

\begin{lemma}\label{BIGGERK}
There exists an integer $k$ such that for all semistable
$X$ sheaves $E$ on $S$ the following holds:

Let $Z$ be a smooth curve in the linear system $|kH|$
such that $E|_Z$ is a vector bundle on $Z$.
If $E|_Z \rarpa{\bar \alpha} \bar M$ is a surjection
onto a $Z$-line bundle $\bar M$ of degree zero,
then $\bar \alpha$ is the restriction of a morphism
$E \rarpa{\alpha} M$ to $Z$, where $M$ is a
$S$-line bundle with $M.H=0$.
\end{lemma}
\proof
This lemma follows from Bogomolov's inequality,
as in the proof of (2) $\Rar$ (3) of \ref{SEMISTABLE}.
For details see \cite{Bog2} theorem 2.3 or \cite{HL} theorem 7.3.5.

\tiptop
{\bf Proof of the positivity theorem \ref{POSITIVE}}\\
Since $\Lcal_B$ is assumed to be the trivial line bundle,
we may pass to a power $\Lcal_B^{\otimes k}$.
Choose $k$ such that lemmas \ref{BIGK} and \ref{BIGGERK} apply.
Now let $F$ be a torsion sheaf supported on a smooth divisor
$Z \in |kH|$ such that the global section
of $\Lcal_B^{\otimes k}$ defined by $F$ is nontrivial.
By construction it is clear that $\Ecal$ restricted to $Z$
is a vector bundle.
Now by theorem \ref{POSCURVE} this yields the
$S$-equivalence of the restriction $\Ecal|_Z$,
and from the proof given here,
it follows that we are in one of the following two cases.

\case{1} All $Z$-vector bundles parameterized by
$B$ are stable and isomorphic.

Let $P$ and $Q$ be two geometric points of the curve $B$.
We consider the following long exact sequence
$$0 \rar \Hom(\Ecal_P, \Ecal_Q(-kH)\ddual)
\rar \Hom(\Ecal_P, \Ecal_Q\ddual) \rarpa{\alpha} $$
$$\rarpa{\alpha} \Hom(\Ecal_P, \Ecal_Q\ddual|_Z)
\rar \Ext^1(\Ecal_P, \Ecal_Q(-kH)\ddual) \rar \, .$$
Since $E_P$ and $E_Q$ are semistable the group
$\Hom(\Ecal_P, \Ecal_Q(-kH)\ddual)$ vanishes.
Hence by lemma \ref{BIGK} the morphism
$\alpha$ is an isomorphism.
The support of the cokernels of the nontrivial morphisms
$\Ecal_P$ to $\Ecal_Q\ddual$ can not change because
they never meet the ample divisor $Z$.
This shows the Jordan-H\"older equivalence of
$\Ecal_Q$ and $\Ecal_P$.

\case{2} All $Z$ vector bundles parameterized by $B$ have a
surjection to a $Z$-line bundle $\bar M$ of degree zero.

First we remark that $\bar M$ is the restriction of an
$S$-line bundle $M$ to $Z$ by lemma \ref{BIGGERK}.
As in the first case, all sheaves parameterized by
$B$ are Jordan-H\"older equivalent. \qed

\section{The Barth morphism}\label{SBAR}
Our construction gives us the moduli space $X$ together with a
finite morphism $X \rarpa{\varphi_1} \pdop^N$,
which we call the Barth morphism.
In this section we show that this morphism for
``surfaces with many lines'' assigns a sheaf $E$
its divisor of jumping curves.
We say that a polarized surface $(S,\Ocal_S(H))$ has many lines
if the linear system $|H|$ is
globally generated and the generic curve of this linear system is
rational.
By adjunction we have $H.(H+K_S) = -2$.
Hence any given rank two sheaf $E$ can be normalised such that
$H.c_1(E) \in \{-3, -2 , -1, 0 \}$ by twisting
$\Ocal_S(H+K_S)^{\otimes k}$.
In the even case the assumptions made in \ref{SEMISTABLE}
are not needed because of the following theorem.
\begin{theorem}
{\bf (Grauert-M\"ulich theorem) }
If $(S,H)$ is a surface with many lines and $E$ a torsion
free rank two sheaf on $S$ with $c_1(E).H$ even,
then $E$ is semistable if and only if the restriction of $E$
to the general curve $l$ in the
linear system $|H|$ is isomorphic to the direct sum of two
isomorphic $l$-line bundles.
\end{theorem}
\proof See \cite{OSS} II  theorem, 2.1.4. or \cite{He2} Lemma 3.6. 

\tiptop
This theorem provides us with a very explicit description of
the dual moduli space $Y$ as we will see in a moment.
We have to distinguish the following two cases:

{\bf Case 1: $H.c_1(E) =-2$}

By the Grauert-M\"ulich theorem, $H^*(E \otimes \Ocal_l)=0$
for $E$ semistable and $l$ general in $|H|$.
Hence we may take $|H|$ to be the dual moduli space.
Recalling our construction we work here with a square root of
the line bundle $\Lcal$. 

{\bf Case 2: $H.c_1(E) =0$ }

In this case we again identify $Y$ with the
complete linear system $|H|$.
Note that for any $\l \in |H|$ the dimension of the
$\Ext^1_{\Ocal_l}(\Ocal_l, \Ocal_l(K_l+H))$ is one.
Hence there is a unique nontrivial extension $\xi_l$.
By assigning to $l$ this torsion sheaf $\xi_l$,
$|H|$ is identified with the dual moduli space $Y$.

\tiptop
Any semistable sheaf $D_E$ defines now a divisor
$D_E$ in the dual moduli space $Y$ consisting of all curves
$l \in |H|$ where $ \otimes \Ocal_l$ is not of the expected type.
Therefore this divisor $D_E$ is called the
{\em divisor of jumping curves}.
By straightforward calculations we find the

\begin{proposition}
The degree $d_E$ of the divisor $D_E$ of jumping curves equals
$$\begin{array}{ll}
d_E=c_2(E)-2+\frac{(K_S-c_1(E)).c_1(E)}{2} & \mbox{ if } H.c_1(E)=-2\,
;\\
\\
d_E=2c_2(E)-2-c_1(E)^2 & \mbox{ if } H.c_1(E)=0\, .\\
\end{array}$$
\end{proposition}

If now we choose $N+1 = { (H^2 +1 ) + d_E \choose d_E }$
curves $\{ l_i \}_{i+0}^N$ in $|H|$ such that
no divisor of degree $d_E$ contains all the $l_i$,
then the duality construction gives us a finite
morphism $X \rarpa{\varphi_1} \pdop^N$.
The name Barth morphism is used because in \cite{Bar}
Barth studied rank-2 vector bundles on the
projective plane via their divisors of jumping lines.
This morphism assigns every semistable rank two sheaf
$E$ its jumping divisor $D_E$.
As a corollary of this construction we have the
\begin{theorem}
The Barth morphism $X \rarpa{\varphi_1} \pdop^N$ is finite.
\end{theorem}

Georg Hein\\
Humboldt-Universit\"at zu Berlin\\
Institut f\"ur Mathematik\\
Burgstr. 26\\
10099 Berlin (Germany) \\

hein@mathematik.hu-berlin.de\\

\end{document}